\documentclass{amsart}
\usepackage{tikz-cd}
\usepackage{latexsym,amssymb}
\usepackage[english]{babel}
\usepackage{bm}
\usepackage{amsfonts, amsthm}
\usepackage{amsmath}
\usepackage{mathrsfs}
\usepackage{multirow}
\usepackage{bm}
\usepackage{pgf,tikz}

\newtheorem{prop}{Proposition}[section]
\newtheorem{lem}[prop]{Lemma}
\newtheorem{cor}[prop]{Corollary}
\newtheorem{thm}[prop]{Theorem}
\newtheorem{conj}[prop]{Conjecture}
\theoremstyle{definition}
\newtheorem{rem}[prop]{Remark}
\newtheorem{defi}[prop]{Definition}
\newtheorem{ex}[prop]{Example}
\newtheorem*{KN}{Crazy Knight's Tour Problem}

\def\Z{\mathbb{Z}}

\numberwithin{equation}{section}

\def\S{\mathbb{S}}
\def\Z{\mathbb{Z}}

\newcommand{\probname}{Crazy Knight's Tour Problem}

\def\G{\Gamma}

\def\H{\mathrm{H}}

\def\S{\mathcal{S}}
\def\E{\mathcal{E}}
\def\D{\mathcal{D}}
\def\B{\mathcal{B}}
\def\F{\mathcal{F}}
\def\R{\mathcal{R}}
\def\C{\mathcal{C}}

\begin{document}

\title[On $\lambda$-fold relative Heffter arrays and biembedding...]{On $\lambda$-fold relative Heffter arrays
and biembedding multigraphs on surfaces}

\author[S. Costa]{Simone Costa}
\address{DICATAM - Sez. Matematica, Universit\`a degli Studi di Brescia, Via
Branze 43, I-25123 Brescia, Italy}
\email{simone.costa@unibs.it}

\author[A. Pasotti]{Anita Pasotti}
\address{DICATAM - Sez. Matematica, Universit\`a degli Studi di Brescia, Via
Branze 43, I-25123 Brescia, Italy}
\email{anita.pasotti@unibs.it}

\begin{abstract}
In this paper we define a new class of partially filled arrays, called $\lambda$-fold relative Heffter arrays, that are a
generalization of the Heffter arrays introduced by Archdeacon in 2015.
After showing the connection of this new concept with several other ones, such as signed magic arrays, graph decompositions and
relative difference families, we determine some necessary conditions and we present existence results for infinite classes of these arrays.
In the last part of the paper we also show that these arrays give rise to biembeddings
of multigraphs into orientable surfaces
and we provide infinite families of such biembeddings.
To conclude, we present a result concerning pairs of $\lambda$-fold relative Heffter arrays and covering surfaces.
\end{abstract}


\maketitle

\section{Introduction}\label{sec:Intro}
An $m \times n$  partially filled (p.f., for short) array on a set $\Omega$ is an $m \times n$ matrix
whose elements belong to $\Omega$ and where we also allow some cells to be empty.
An interesting class of p.f. arrays, called Heffter arrays, has been introduced by Dan Archdeacon in \cite{A}
and then generalised in \cite{RelH} as follows.

\begin{defi}\cite{RelH}\label{def:RelativeH}
Let $v=2nk+t$ be a positive integer, where $t$ divides $2nk$, and let $J$ be the subgroup of $\Z_v$ of order $t$.
 A $\H_t(m,n; s,k)$ \emph{Heffter array  over $\Z_v$ relative to $J$} is an $m\times n$ p.f.  array
 with elements in $\Z_v$ such that:
\begin{itemize}
\item[($\rm{a})$] each row contains $s$ filled cells and each column contains $k$ filled cells;
\item[($\rm{b})$] for every $x\in \Z_{v}\setminus J$, either $x$ or $-x$ appears in the array;
\item[($\rm{c})$] the elements in every row and column sum to $0$ (in $\Z_v$).
\end{itemize}
\end{defi}
The classical concept of Heffter array
introduced by Archdeacon in \cite{A} corresponds to a relative Heffter array with $t=1$,
namely when $J$ is the trivial subgroup of $\Z_{2nk+1}$.
In general $t$ is omitted when its value is $1$, hence a classical Heffter array is denoted by $\H(m,n;s,k)$.
If we have a square array, that is if $m=n$, then $s=k$ and a  $ \H_t(n,n; k,k)$ will be simply denoted by  $\H_t(n; k)$.
 A relative Heffter array is called \emph{integer} if Condition ($\rm{c}$) in Definition \ref{def:RelativeH} is
strengthened so that the elements in every row and in every
column, viewed as integers in
$\pm\left\{ 1, \ldots, \left\lfloor \frac{2nk+t}{2}\right\rfloor \right\} ,$
sum to zero in $\Z$.

Classical Heffter arrays and their generalization given in Definition \ref{def:RelativeH} are considered interesting and worthy of study
also because they have several applications.
Among these we recall a construction of complete graph embeddings, a topic studied also in view of its connection with the ``Heawood Map Colouring Conjecture'',
see \cite{GG} and \cite{R} for more details on this relation and \cite{CHK} for other constructions of complete graph embeddings.
So, there are some recent papers in which Heffter arrays are investigated to obtain new face 2-colorable embeddings, briefly biembeddings
(see \cite{A,CDDYbiem, CMPPHeffter, CPPBiembeddings, DM}). On the other hand there are also papers entirely dedicated to the existence problem
(see \cite{ABD, ADDY,  BCDY, CDDY, RelH, DW, MP}). In particular, in \cite{ADDY,DW} the authors verify the existence of a square integer Heffter array for all admissible orders,
proving the following theorem.
\begin{thm}\label{thm:squareinteger}
  There exists an integer $\H(n;k)$ if and only if $3\leq k\leq n$ and $nk\equiv 0,3\pmod 4$.
\end{thm}
Also, the existence problem of \emph{relative} integer Heffter arrays has been investigated in \cite{RelH, CPPBiembeddings, MP}.
The main results are the following.
\begin{thm}\label{thm:esistenzaRel}\cite{RelH}
Let $3\leq k\leq n$ with $k\neq 5$.
  There exists an integer $\H_k(n;k)$ if and only if one of the following holds:
  \begin{itemize}
    \item[(1)] $k$ is odd and $n\equiv 0,3\pmod 4$;
    \item[(2)] $k\equiv 2\pmod 4$ and $n$ is even;
    \item[(3)] $k\equiv 0\pmod 4$.
  \end{itemize}
Furthermore, there exists an integer $\H_5(n;5)$ if
$n\equiv 3\pmod 4$ and it does not exist if $n\equiv 1,2\pmod 4$.
\end{thm}

\begin{prop}\label{n,2n}\cite{CPPBiembeddings}
For every odd $n\geq 3$ there exists an integer $\H_n(n;3)$ and an integer $\H_{2n}(n;3)$.
\end{prop}

\begin{thm}\label{main}\cite{MP}
Let $m,n,s,k$ be integers such that $4\leq s\leq n$, $4\leq k\leq m$ and $ms=nk$. Let $t$ be a divisor of $2nk$.
\begin{itemize}
\item[(1)] If $s\equiv k \equiv 0 \pmod 4$, then  there exists an integer $\H_t(m,n;s,k)$.
\item[(2)] If $s \equiv 2 \pmod 4$ and $k\equiv 0\pmod 4$, then an integer $\H_t(m,n;s,k)$ exists if and only if $m$ is even.
\item[(3)]  If $s \equiv 0 \pmod 4$ and $k\equiv 2\pmod 4$, then an integer $\H_t(m,n;s,k)$ exists if and only if $n$ is even.
\item[(4)] Suppose that $m$ and $n$ are both even. If $s\equiv k \equiv 2 \pmod 4$, then  there exists an integer $\H_t(m,n;s,k)$.
\end{itemize}
\end{thm}

In this paper we propose a further natural generalization of Heffter arrays, which is related to signed magic arrays, difference families,
graphs decompositions and biembeddings, as shown in the following sections.

\begin{defi}\label{def:lambdaRelative}
Let $v=\frac{2nk}{\lambda}+t$ be a positive integer,
where $t$ divides $\frac{2nk}{\lambda}$,  and
let $J$ be the subgroup of $\Z_{v}$ of order $t$.
 A  $\lambda$-\emph{fold Heffter array $A$ over $\Z_{v}$ relative to $J$}, denoted by $^\lambda\H_t(m,n; s,k)$, is an $m\times n$ p.f.  array
 with elements in $\Z_{v}$ such that:
\begin{itemize}
\item[($\rm{a_1})$] each row contains $s$ filled cells and each column contains $k$ filled cells;
\item[($\rm{b_1})$] the multiset $\{\pm x \mid x \in A\}$ contains $\lambda$ times each element of $\Z_v\setminus J$;
\item[($\rm{c_1})$] the elements in every row and column sum to $0$ (in $\Z_v$).
\end{itemize}
\end{defi}
Trivial necessary conditions for the existence of a $^\lambda \H_t(m,n; s,k)$ are $ms=nk$, $2\leq s \leq n$ and $2\leq k \leq m$.
It is easy to see that condition $(\rm{b_1})$ of Definition \ref{def:lambdaRelative} asks that
for every $x\in \Z_{v}\setminus J$ with $x$ different from the involution the sum of the occurrences of $x$ and $-x$ in the array is $\lambda$,
while if the involution exists and it does not belong to $J$ then it has to appear exactly $\frac{\lambda}{2}$ times,
also no element of $J$ appears in the array.

If $^\lambda \H_t(m,n; s,k)$ is a square array, then it will be denoted by $^\lambda \H_t(n;k)$.
Note that if $\lambda=1$ we get the concept of a relative Heffter array given in Definition \ref{def:RelativeH}, also if $\lambda=t=1$
 we find again the classical concept of Heffter array introduced in \cite{A}.
 Note also that if $\lambda=1$, then $3\leq s \leq n$ and $3\leq k \leq m$.

\begin{ex}\label{ex:1}
  Below we have a $^3\H_2(4;3)$ whose elements belong to $\Z_{10}$ and a $^4\H_4(4;2)$ whose elements belong to $\Z_8$.
  \center{
  $\begin{array}{|r|r|r|r|}
\hline  & 1 & 2  & -3  \\
\hline 4  &  & 4 & 2  \\
\hline -3  & 2 &  &  1 \\
\hline   -1 & -3 & 4 &  \\
\hline
\end{array}$
\quad\quad\quad\quad
  $\begin{array}{|r|r|r|r|}
\hline  1 & -1 &   &   \\
\hline -1  & 1 &  &   \\
\hline   &  & 3 &  -3 \\
\hline   &  & -3 & 3 \\
\hline
\end{array}$
}
\end{ex}
A $\lambda$-fold relative Heffter array is called \emph{integer} if condition ($\rm{c_1}$) in Definition \ref{def:lambdaRelative}
is strengthened so that the elements in every row and every column, seen as integer in $\pm\{1,2,\ldots, \lfloor\frac{v}{2}\rfloor\}$, sum to $0$ in $\Z$.
Note that the $^3\H_2(4;3)$ of Example \ref{ex:1} is not an integer Heffter array, while the $^4\H_4(4;2)$ in the same example is integer.
\begin{ex}\label{ex:inv}
Below we have an integer $^2\H_1(5;3)$. Here we are working in $\Z_{16}$ and the involution does not belong to the subgroup
of order $t=1$, that is the trivial subgroup.
Hence $8$ appears exactly once in the array.

  \center{
  $\begin{array}{|r|r|r|r|r|}
\hline 3 & 5 & -8  & &  \\
\hline & 2 & 5 & -7 &   \\
\hline & & 3  & 1 &  -4 \\
\hline  -4 &  &  & 6 & -2 \\
\hline   1 & -7 & &  & 6 \\
\hline
\end{array}$

}
\end{ex}

\begin{ex}\label{ex:2}
  Below we have an integer $^2\H(6;4)$ whose elements belong to $\Z_{25}$.
  \center{
  $\begin{array}{|r|r|r|r|r|r|}
\hline & & 1 & -2  & -5 & 6  \\
\hline & & -3  & 4 & 7 & -8  \\
\hline -9 & 10  &  &  &  1 & -2 \\
\hline   11 & -12 & & &  -3 & 4   \\
\hline 5 & -6 &-9 & 10  &  &   \\
\hline  -7 & 8 & 11 & -12 & &   \\
\hline
\end{array}$}
\end{ex}

In Section \ref{sec:magic} we show how signed magic arrays are nothing but a very particular case
of $2$-fold Heffter arrays, while in Section \ref{sec:DF} we present the connections with cyclic
cycle decompositions of the complete multipartite multigraph, this connection is obtained via
relative difference families, that are also related to  $\lambda$-fold relative Heffter arrays.
In Section \ref{sec:NC} after having determined some necessary conditions for the existence of such arrays
and a non existence result, we prove Theorem \ref{thm:proiezione} which will turn out to be very useful to obtain
new infinite classes of $\lambda$-fold relative Heffter arrays starting from know results about (relative) Heffter arrays.

The last part of the paper is dedicated to biembeddings. In Section 5, we show that also this generalization of Heffter array, as well as
the classical concept introduced by Archdeacon, is useful for finding biembeddings of cycle
decompositions. In particular, we prove how a $\lambda$-fold relative
Heffter array gives rise to a cellular biembedding of a pair of cyclic cycle decompositions
of a complete multipartite multigraph into an orientable surface and we present results
about these biembeddings.
To conclude,  using the Heffter arrays of Theorem 4.4 and considering the related embeddings
of multigraphs into two orientable surfaces $\Sigma$ and $\Sigma'$, we show that $\Sigma'$ is a covering space
of $\Sigma$.

\section{Relation with signed magic arrays}\label{sec:magic}
We have to point out that in the particular case in which $\lambda=2$
and $t=1$, there is a strong relation with \emph{signed magic arrays}, introduced in \cite{KSW} and
also called signed magic rectangles.
\begin{defi}\label{def:SMR}
A \emph{signed magic array} $SMA(m,n;s,k)$ is an $m\times n$ array with entries from $X$,
where $X=\left\{0,\pm1,\pm2,\ldots, \pm \frac{nk-1}{2}\right\}$ if $nk$ is odd and $X=$\\ $\left\{\pm1,\pm2,\ldots,\pm \frac{nk}{2}\right\}$
if $nk$ is even, such that
\begin{itemize}
\item[($\rm{a_2})$] each row contains $s$ filled cells and each column contains $k$ filled cells;
\item[($\rm{b_2})$] every integer from the set $X$ appears exactly once in the array;
\item[($\rm{c_2})$] the sum of each row and of each column is zero.
\end{itemize}
\end{defi}
A $SMA(n,n;k,k)$ is also called  a \emph{signed magic square} and it is denoted by $SMA(n;k)$.
Looking at Definitions \ref{def:lambdaRelative} and \ref{def:SMR} it is easy to see that there is the following relation
between $2$-fold Heffter arrays and signed magic arrays.

\begin{rem}\label{rem:magicHeffter}
A signed magic array $SMA(m,n;s,k)$ with $nk$ even is an integer $^2\H(m,n;s,k)$.
We point out that in general the converse is not true. For instance, the $^2\H(6;4)$ of Example \ref{ex:2} is not a signed magic square.
On the other hand,  in the particular case in which either $s=2$ or $k=2$ an
integer $^2\H(m,n;s,k)$ is nothing but a $SMA(m,n;s,k)$. In fact, under this hypothesis,
 each row (respectively each column) contains two elements of the form $x$ and $-x$.
Hence each  element of $X=\left\{\pm1,\pm2,\ldots,\pm \frac{nk}{2}\right\}$
appears exactly once in the array, thus condition $(\rm{b_2})$ of Definition \ref{def:SMR} is satisfied.
Since the Heffter array is integer also condition $(\rm{c_2})$ is satisfied, and  condition
$(\rm{a_2})$ trivially holds.
\end{rem}

Hence before investigating the existence of $2$-fold Heffter arrays it is natural to look for known results about
signed magic arrays $SMA(m,n;s,k)$ with $nk$ even.
As far as we know, the following are all known results.

\begin{thm}\cite{KSW}\label{thm:SMS}
  An $SMA(m,n;n,m)$ exists if and only if one of the following holds:
  \begin{itemize}
   \item[(1)]$m=n=1$;
   \item[(2)] $m=2$ and $n\equiv 0,3\pmod 4$;
   \item[(3)] $n=2$ and $m\equiv 0,3\pmod 4$;
   \item[(4)] $m,n>2$.
\end{itemize}
\end{thm}

About the rectangular case with empty cells there are  the following results.

\begin{thm}\cite{KE}\label{thm:magic2}
  There exists an $SMA(m,n;s,2)$ if and only if either $m=2$ and $n=s\equiv 0,3 \pmod 4$
  or $m,s\geq 3$ and $ms=2n$.
\end{thm}

\begin{thm}\cite{KLE}\label{thm:magic3}
  There exists an $SMA(m,n;s,3)$ if and only if $3\leq m,s\leq n$ and $ms=3n$.
\end{thm}

\begin{thm}\cite{MP2}\label{thm:magicMP}
Let $s,k$ be even integers with $s,k \geq 4$.
There exists an $SMA(m,n;s,k)$ if and only if $4\leq s\leq n$, $4\leq k \leq m$ and $ms=nk$.
\end{thm}

In \cite{KSW} the authors completely solve the square case as shown in the following result.
\begin{thm}\cite{KSW}\label{SMS}
There exists an $SMA(n;k)$ if and only if $n=k=1$ or $3\leq k\leq n$.
\end{thm}

In Section \ref{sec:NC} we will use some of these results to
 give a complete solution to the existence problem of a $2$-fold (integer) Heffter array
 $A$ in each of the following cases: $A$ is square (see Theorem \ref{lambda2}),
 each row and each column of $A$ contains an even number of filled cells
 (see Theorems \ref{thm:skeveninteger} and \ref{thm:skeven}),
 $A$ has no empty cells (see Theorem \ref{thm:noemptycell}).


\section{Relations with graph decompositions and difference families}\label{sec:DF}
First we recall some definitions and we give some notation.
Given a graph $\G$, we denote by $V(\G)$ and $E(\G)$ the vertex-set and the edge-set of $\G$, respectively,
and by $^\lambda \G$ the multigraph obtained from $\G$ by repeating each edge $\lambda$ times.
Also by $K_v$, $K_{q \times r}$ and $C_k$ we represent the complete graph on $v$ vertices,
the complete multipartite graph with $q$ parts each of size $r$ and the cycle of length $k$, respectively.
Given a subgraph $\Gamma$ of a graph $K$, a $\G$-decomposition of $K$ is a set of graphs, all isomorphic to $\G$,
whose edge-sets partition the edge-set of $K$, see for instance \cite{BEZ}.
Given an additive group $G$, a $\G$-decomposition $\D$ of a graph $K$ is $G$-regular if, up to isomorphisms,
$V(K)=G$ and for any $B\in {\D}$ also the graph $B+g$ is a block of ${\D}$ for any $g\in G$.
Here we are interested in cyclic cycle decompositions, namely $\G$-decompositions regular under the cyclic group
with $\G$ a cycle.
Details about regular cycle decompositions can be found in \cite{B04}, where it is shown that difference families are
a very useful tool for finding such decompositions.
We  recall the definition,  see also \cite{AB,B98}.
\begin{defi}
Let $\G$ be a graph with vertices in an additive group $G$. The multiset
$\Delta \G=\{\pm(x-y) \mid \{x,y\}\in E(\G)\}$
is called the \emph{list of differences} from $\G$.
\end{defi}
More generally, given a set $\mathcal{W}$ of graphs with vertices in $G$, by $\Delta\mathcal{W}$ one
means the union (counting multiplicities) of all multisets $\Delta\G$, where $\G\in \mathcal{W}$.

\begin{defi}\label{def:DF}
  Let $J$ be a subgroup of an additive group $G$ and let $\G$ be a graph.
  A collection $\F$ of graphs isomorphic to $\G$ and with vertices in $G$
  is said to be a $(G,J,\G,\lambda)$-\emph{difference family} (briefly, DF) \emph{over $G$ relative to} $J$
if each element of $G\setminus J$ appears exactly $\lambda$ times in the list of differences of $\F$
while no element of $J$ appears there.
\end{defi}

If $t$ is a divisor of $v$, a $(\Z_v,\frac{v}{t}\Z_v,\G,\lambda)$-DF, where $\frac{v}{t}\Z_v$ denotes the subgroup of $\Z_v$
of order $t$, is simply denoted by $(v,t,\G,\lambda)$-DF.
The connection between relative difference families and decompositions of a complete multipartite multigraph
is given by the following result.

\begin{prop}\label{thm:basecycles}\cite[Proposition 2.6]{BP}
If $\F=\{B_1,\ldots,B_\ell\}$ is a $(G,J,\G,\lambda)$-DF, then $\B=\{B_i+g \mid i=1,\ldots,\ell; g\in G\}$
is a $G$-regular $\G$-decomposition of $^\lambda K_{q\times r}$, where $q=|G:J|$ and $r=|J|$.
\end{prop}

Now, in order to present the connection between $\lambda$-fold relative Heffter arrays and relative difference families,
we have to introduce the concept of \emph{simple ordering}, which also plays a fundamental role in the study of the relation
with biembeddings, as shown in Section \ref{sec:biembedding}.

Henceforward, given two integers $a\leq b$, we denote by $[a,b]$ the interval containing the integers $a,a+1,\ldots,b$.
If $a>b$, then $[a,b]$ is empty.

If $A$ is an  $m\times n$ p.f.  array we define the skeleton of $A$, denoted by $skel(A)$, to be the set of the filled positions of $A$;
we will also denote by $\E(A)$ the multiset of the elements of $skel(A)$.
Then, we name by $\overline{R}_1,\ldots,\overline{R}_m$ and by $\overline{C}_1,\ldots,\overline{C}_n$, respectively,
 the rows and the columns of $A$. Analogously, by $\E(\overline{R}_i)$ and $\E(\overline{C}_j)$ we mean the multisets of elements
 of the $i$-th row and of the $j$-th column, respectively, of $A$.
We consider  the elements of $\E(A)$  indexed by the set $skel(A)$ and the elements of $\E(\overline{R}_i)$
 (resp. $\E(\overline{C}_j)$) indexed by the set $skel(A)\cap \overline{R}_i$ (resp. $skel(A)\cap \overline{C}_j$).
For example let $A$  be the $^4\H_4(4;2)$ constructed in Example \ref{ex:1}.
Here $\E(A)=\{-1,-1,1,1,-3,-3,3,3\}$, we can view this multiset as the set
$\{t_b \mid b\in skel(A)\}=\{t_{(1,1)},t_{(1,2)},t_{(2,1)},t_{(2,2)},t_{(3,3)},t_{(3,4)},t_{(4,3)},t_{(4,4)}\}$,
where $t_{(1,1)}=1$, $t_{(1,2)}=-1$, $t_{(2,1)}=-1$, $t_{(2,2)}=1$, $t_{(3,3)}=3$, $t_{(3,4)}=-3$,
$t_{(4,3)}=-3$ and $t_{(4,4)}=3$.
Given a finite multiset $T=[t_{b_1},\dots,t_{b_k}]$ whose (not necessarily distinct) elements are indexed by a set $B=\{b_1,\dots,b_k\}$ and a cyclic permutation $\alpha$ of $B$ we say that the list $\omega_{\alpha}=(t_{\alpha(b_1)},t_{\alpha(b_2)},\ldots,t_{\alpha(b_k)})$ is the ordering of the elements of $T$ associated to $\alpha$. In case the elements of $T$ belong to an abelian group $G$, we define $s_i=\sum_{j=1}^i t_{\alpha(b_j)}$, for any $i\in[1,k]$, to be the $i$-th partial sum of $\omega_{\alpha}$ and we set $\S(\omega_{\alpha})=(s_1,\ldots,s_k)$.
The ordering $\omega_{\alpha}$ is said to be \emph{simple} if $s_b\neq s_c$ for all $1\leq b <  c\leq k$ or,
equivalently, if there is no proper subsequence of $\omega_{\alpha}$ that sums to $0$.
Note that if $\omega_{\alpha}$ is a simple ordering so is $\omega_{\alpha^{-1}}=(t_{\alpha(b_k)},t_{\alpha(b_{k-1})},\ldots,t_{\alpha(b_1)})$. In the following, if there are no ambiguities, we will just write $\omega$ for $\omega_{\alpha}$ and $\omega^{-1}$ for $\omega_{\alpha^{-1}}$ omitting the dependence on $\alpha$.
We point out that there are several interesting problems and conjectures about distinct partial sums.
For instance several years ago Alspach made the following conjecture, whose validity would shorten some cases of known proofs about the existence of cycle decompositions.

\begin{conj}\label{Conj:als}
Let $A\subseteq \mathbb{Z}_v\setminus\{0\}$ such that
$\sum_{a\in A}a\neq0$. Then there exists an ordering of the elements of $A$ such that the partial sums are all
distinct and nonzero.
\end{conj}
Results on this conjecture have been obtained in \cite{A,BH, CP,HOS}.
Other problems related to the previous one have been proposed in
 \cite{ AL, ADMS, CMPPSums, O}.
Given an $m \times n$ p.f. array $A$, by $\alpha_{\overline{R}_i}$ we will denote a cyclic permutation of the nonempty cells of $\overline{R}
_i$ and by $\omega_{\overline{R}_i}$ (omitting the dependence on $\alpha_{\overline{R}_i}$) the associated ordering of $\E(\overline{R}_i)$. Similarly we define $\alpha_{\overline{C}_j}$ and $\omega_{\overline{C}_j}$.

 If for any $i\in[1, m]$ and for any $j\in[1,n]$, the orderings $\omega_{\overline{R}_i}$ and $\omega_{\overline{C}_j}$ are simple, we define the permutation $\alpha_r$ of $skel(A)$
 by $\alpha_{\overline{R}_1}\circ \ldots \circ\alpha_{\overline{R}_m}$ and we say that the associated action $\omega_r$ on $\E(A)$ is the simple ordering for the rows.
 Similarly we define $\alpha_c=\alpha_{\overline{C}_1}\circ \ldots \circ\alpha_{\overline{C}_n}$ and the simple ordering $\omega_c$ for the columns.
Moreover, by \emph{natural ordering} of a row (column) of $A$ we mean the ordering from left to right (from top to
bottom).
A p.f. array $A$ on an abelian group $G$ is said to be
\begin{itemize}
\item \emph{simple} if each row and each column of $A$ admits a simple ordering;
\item \emph{globally simple} if the natural ordering of each row and each column of $A$ is simple.
\end{itemize}

\begin{rem}\label{rem:semplice}
If $s$ and $k$ do not exceed $3$, then every ordering for the rows and the columns of a  $^\lambda \H_t(m,n;s,k) $ is simple.
\end{rem}

Now we are ready to explain the relation between $\lambda$-fold relative Heffter arrays and relative difference families.
We recall that by $C_k$ we denote the cycle on $k$ vertices.

\begin{prop}\label{from Heffter to DF}
  If there exists a simple $^\lambda \H_t(m,n;s,k) $, then there exists a\\ $\left(\frac{2ms}{\lambda}+t,t,C_s,\lambda\right)$-DF and a
  $\left(\frac{2nk}{\lambda}+t,t,C_k,\lambda\right)$-DF.
\end{prop}
\begin{proof}
Let $A$ be a simple $^\lambda \H_t(m,n;s,k) $. Then, for any $i\in [1,m]$, there exists a simple ordering $\omega_{\overline{R}_i}$
of the $i$-th row of $A$. Hence, from each row of $A$ we can construct an $s$-cycle whose vertices in $\Z_{\frac{2ms}{\lambda}+t}$
are the partial sums of $\omega_{\overline{R}_i}$. Denoting by $\F_s$ the set of $m$ $s$-cycles so constructed starting from the rows of $A$,
we obtain $\Delta{\F}_s=\pm{\E(A)}$. On the other hand, since $A$ is a $^\lambda \H_t(m,n;s,k) $,
by condition ($\rm{b_1}$) of Definition \ref{def:lambdaRelative}, we have that $\pm{\E(A)}$
contains exactly $\lambda$ times each element of $\Z_{\frac{2ms}{\lambda}+t}\setminus J$,
where $J$ is the subgroup of order $t$ of $\Z_{\frac{2ms}{\lambda}+t}$, and no element
of $J$ is contained in  $\pm{\E(A)}$. Thus, $\F_s$ is a $\left(\frac{2ms}{\lambda}+t,t,C_s,\lambda\right)$-DF.

In a similar way, starting from the columns of $A$, one can construct a\\ $\left(\frac{2nk}{\lambda}+t,t,C_k,\lambda\right)$-DF.
\end{proof}

\begin{ex}
  Starting from the array $A=$ $^3\H_2(4;3)$ given in Example \ref{ex:1} we can construct two
   $(10,2,C_3,3)$-DFs. Since $k=3$ by Remark \ref{rem:semplice} every ordering is simple. So considering for instance the natural
   ordering of each row and each column of $A$ we obtain the following $3$-cycles:
   \begin{align*}
C^{\bar{R}_1}&=(1,3,0), & C^{\bar{C}_1}&=(4,1,0),\\
C^{\bar{R}_2} &=(4,8,0), & C^{\bar{C}_2}&=(1,3,0), \\
 C^{\bar{R}_3} &=(-3,-1,0), & C^{\bar{C}_3}&=(2,6,0), \\
 C^{\bar{R}_4} &=(-1,-4,0), & C^{\bar{C}_4}&=(-3,-1,0).
\end{align*}
Set $\F^{\bar{R}}_3=\{C^{\bar{R}_i}\mid i\in[1,4]\}$ and $\F^{\bar{C}}_3=\{C^{\bar{C}_i}\mid i\in[1,4]\}$;
by the construction of the cycles it immediately follows that
$\Delta \F^{\bar{R}}_3$ and $\Delta \F^{\bar{C}}_3$
contain exactly 3 times each element of $\Z_{10}\setminus 5\Z_{10}$
and no element of $5\Z_{10}$.
Hence $\F^{\bar{R}}_3$ and $\F^{\bar{C}}_3$ are two $(10,2,C_3,3)$-DFs.
\end{ex}

To conclude this section we present the connection between $\lambda$-fold relative Heffter arrays
and cyclic cycle decompositions that follows by previous results.
\begin{prop}\label{HeffterToDecompositions}
  Let $A$ be a simple $^\lambda \H_t(m,n;s,k)$ with respect to the simple orderings $\omega_r$ and
  $\omega_c$. Then:
  \begin{itemize}
    \item[(1)] there exists a cyclic $s$-cycle decomposition $\D_{\omega_r}$ of  $^\lambda K_{(\frac{2ms}{\lambda t}+1)\times t}$;
    \item[(2)] there exists a cyclic $k$-cycle decomposition $\D_{\omega_c}$ of  $^\lambda K_{(\frac{2nk}{\lambda t}+1)\times t}$.
  \end{itemize}
\end{prop}

\begin{proof}
(1) and (2) follow from Propositions \ref{thm:basecycles} and \ref{from Heffter to DF}.
\end{proof}

\section{Necessary conditions and  existence results for $^{\lambda} \H_t(m,n;s,k)$}\label{sec:NC}
In this section we firstly determine some necessary conditions for the existence of a $\lambda$-fold relative Heffter array.
Then we  present a result that allows us to obtain infinite classes of
$\lambda$-fold relative Heffter arrays starting from know results on (relative) Heffter arrays.
To conclude we present also a general recursive construction.

We start by recalling that, by definition, $t$ divides $\frac{2nk}{\lambda}$. Also,
conditions ($\rm{a_1})$ and ($\rm{b_1})$ of Definition \ref{def:lambdaRelative} imply that
if there exists a $^{\lambda} \H_t(m,n;s,k)$ with either $s=2$ or $k=2$, then $\lambda$ has to be even.

\begin{prop}\label{prop:necctrivial}
Suppose that there exists a $^\lambda \H_t(m,n;s,k)$ and set $v=\frac{2nk}{\lambda}+t$.
 If either $v$ is odd or $v$ and $t$ are even  then $\lambda$ has to be a divisor of $nk$.
\end{prop}
\begin{proof}
Let $A$ be a $^\lambda \H_t(m,n;s,k)$.
Clearly if $v$ is odd, $\Z_v$ does not have the involution.
Also, if $v$ and $t$ are even, the involution of $\Z_v$ belongs to the subgroup $J$ of order $t$ of $\Z_v$.
Hence in both cases, the array does not contain the involution.
So for every  $x \in \Z_v\setminus J$, the total number of occurrences of $x$ and $-x$ in $A$ is $\lambda$.
The thesis follows.
\end{proof}

\begin{prop}\label{prop:nonexistence}
If $\lambda\equiv 2 \pmod 4$, $v=\frac{2nk}{\lambda}+t\equiv2\pmod 4$ and $t$ is odd, then a  $^{\lambda} \H_t(m,n;s,k)$
cannot exist.
\end{prop}
\begin{proof}
By contradiction, we suppose that there exists a $^{\lambda} \H_t(m,n;s,k)$, say $A$.
Since by hypothesis $v$ is even, in order for each row to sum to zero modulo $v$, each row must contain an even number of odd numbers.
In particular the entire array $A$ contains an even number of odd numbers.
Note that, since $v\equiv 2\pmod 4$, the involution $i$ of $\Z_v$ is an odd integer.
Also, since $t$ is odd, $i$ does not belong to the subgroup of order $t$ of $\Z_v$,
hence $i$ appears in the array.
So if $x$ is an odd element different from $i$ then the sum of occurrences  of $x$ and $-x$
in $A$ is $\lambda$, while $i$ appears $\frac{\lambda}{2}$ times in $A$.
Since $\lambda\equiv2\pmod 4$, the array $A$ contains an odd number of odd numbers, that is a contradiction.
\end{proof}

If we focus on the \emph{integer} case,
reasoning as in the proof of Proposition 3.1 of \cite{RelH}, one can obtain also the following result.

\begin{prop}\label{prop:necc}
Suppose that there exists an integer $^\lambda \H_t(m,n;s,k)$ with $\lambda$ odd.
\begin{itemize}
\item[(1)] If $t$ divides $\frac{nk}{\lambda}$, then
$$\frac{nk}{\lambda}\equiv 0 \pmod 4 \quad \textrm{ or } \quad \frac{nk}{\lambda}\equiv -t \equiv \pm 1\pmod 4.$$
\item[(2)] If $t=\frac{2nk}{\lambda}$, then $s$ and $k$ must be even.
\item[(3)] If $t\neq \frac{2nk}{\lambda}$ does not divide $\frac{nk}{\lambda}$, then
$$\frac{2nk}{\lambda}+t\equiv 0 \pmod 8.$$
\end{itemize}
\end{prop}

The following theorem is a very useful tool for obtaining new $\lambda$-fold
relative Heffter arrays starting from known existence results.

\begin{thm}\label{thm:proiezione}
If there exists an $^\alpha \H_t(m,n;s,k)$ then there exists a $^{\lambda\alpha} \H_\frac{t}{\lambda}(m,n;s,k)$ for any divisor $\lambda$
of $t$.
\end{thm}
\begin{proof}
Let $A=$ $^\alpha\H_t(m,n;s,k)$, clearly $t$ divides $\frac{2nk}{\alpha}$, and set $v=\frac{2nk}{\alpha}+t$.
Let $\lambda$ be a divisor of $t$, and hence of $\frac{2nk}{\alpha}$.
Set $B$ to be the $m\times n$ array obtained from $A$ by considering each element $x$ contained in $A$
as an element of $\Z_\frac{v}{\lambda}$.
We want to show that $B$ is a  $^{\lambda\alpha} \H_\frac{t}{\lambda}(m,n;s,k)$.
Note that $\frac{t}{\lambda}$ divides $\frac{2nk}{\lambda\alpha}$. Now we prove that $B$ satisfies the conditions
of Definition \ref{def:lambdaRelative}.
It is trivial that  condition ($\rm{a_1}$) is satisfied.
It is easy to see that $B$ satisfies also condition ($\rm{c_1}$), since the elements of every row and column sum $0$
in $\Z_v$ and hence also in $\Z_\frac{v}{\lambda}$.
So we are left to consider condition ($\rm{b_1}$), namely we have to prove that the multiset
$\{\pm x \mid x\in B\}$ contains $\lambda\alpha$ times each element of $\Z_\frac{v}{\lambda}\setminus J$,
where $J$ is the subgroup of $\Z_\frac{v}{\lambda}$ of order $\frac{t}{\lambda}$.
Note that we can write $\pm \E(A)$ as the following disjoint union
$$\pm \E(A)=\alpha X_1\ \dot\cup\ \alpha X_2\ \dot\cup\  \ldots \  \dot\cup \ \alpha X_\frac{2nk}{\lambda\alpha}$$
where $X_i=\{x_i,x_i+\frac{v}{\lambda},x_i+2\frac{v}{\lambda},\ldots,x_i+(\lambda-1)\frac{v}{\lambda}  \}$, for $i=1,\ldots, \frac{2nk}{\lambda\alpha}$,
that is $X_i$ contains all the elements of $\Z_v$
equivalent to $x_i\pmod{\frac{v}{\lambda}}$
and where by $\alpha X_i$ we mean the multiset containing $\alpha$ times each element of $X_i$.
It follows that $\pm \E(B)$ contains $\lambda\alpha$ times each element $x_i$ of $\Z_\frac{v}{\lambda}\setminus J$.
Clearly since $x_i\in \pm \E(A)$, then $x_i$ does not belong to the subgroup of order $t$ of
$\Z_v$, hence $x_i$, considered as element of $\Z_\frac{v}{\lambda}$, does not belong to $J$.
This implies that since no element of the subgroup of order $t$ of $\Z_v$ appears in $A$
then no element of the subgroup $J$ of $\Z_\frac{v}{\lambda}$ appears in $B$.
The thesis follows.

\end{proof}

In this section we will apply the following immediate consequence.
\begin{cor}\label{cor:proiezione}
If there exists a $\H_t(m,n;s,k)$ then there exists a $^\lambda \H_\frac{t}{\lambda}(m,n;s,k)$ for any divisor $\lambda$
of $t$.
\end{cor}

\begin{rem}
  Note that given an integer $^\alpha\H_t(m,n;s,k)$,
this does not imply that the $^{\lambda\alpha} \H_\frac{t}{\lambda}(m,n;s,k)$
obtained following the proof of Theorem \ref{thm:proiezione} is integer too,
as shown in the following example.
\end{rem}

\begin{ex}\label{ex:proiezione}
Here we present an application of Theorem \ref{thm:proiezione}. Consider the following integer $\H_2(13;3)$ whose elements
belong to
$\Z_{80}$:
{\small{
\begin{center}
  $\begin{array}{|r|r|r|r|r|r|r|r|r|r|r|r|r|}
\hline -12 & 26 &   & & & & & & & & & & -14  \\
\hline -27 & -11 &  38 & & & & & & & & & &  \\
\hline  & -15 & -10  & 25 & & & & & & & & &  \\
\hline  &  & -28  & -9 & 37 & & & & & & & &  \\
\hline  &  &   & -16 & -8 & 24 & & & & & & &  \\
\hline  &  &   & & -29 & -7 & 36 & & & & & &  \\
\hline  &  &   & & & -17 & -13 & 30 & & & & &  \\
\hline  &  &   & & & & -23 & 5 & 18 & & & &   \\
\hline  &  &   & & & & & -35 & 4 & 31 & & &   \\
\hline  &  &   & & & & & & -22 & 3 & 19 & &   \\
\hline  &  &   & & & & & & & -34 & 2 & 32 &   \\
\hline  & &   & & & & & & & & -21 & 1 & 20  \\
\hline 39 & &   & & & & & & & & & -33 & -6  \\
\hline
\end{array}$
\end{center}}}
If each element of $\H_2(13;3)$ is now considered as an element of $\Z_{40}$
the same array is a $^2\H(13;3)$. In order to see that this new array is not integer,
for each $x$ in $^2\H(13;3)$ we have to consider the element $y\equiv x \pmod{40}$ with $y\in [-20,20]$.
In this way we obtain the following array, say $B$:
\begin{center}
  $\begin{array}{|r|r|r|r|r|r|r|r|r|r|r|r|r|}
\hline -12 & -14 &   & & & & & & & & & & -14  \\
\hline 13 & -11 &  -2 & & & & & & & & & &  \\
\hline  & -15 &   -10  & -15 & & & & & & & & &  \\
\hline  &  & 12 & -9 & -3 & & & & & & & &  \\
\hline  &  &   & -16 & -8 & -16 & & & & & & &  \\
\hline  &  &   & & 11 & -7 & -4 & & & & & &  \\
\hline  &  &   & & & -17 & -13 & -10 & & & & &  \\
\hline  &  &   & & & & 17 & 5 & 18 & & & &   \\
\hline  &  &   & & & & & 5 & 4 & -9 & & &   \\
\hline  &  &   & & & & & & 18 & 3 & 19 & &   \\
\hline  &  &   & & & & & & & 6 & 2 & -8 &   \\
\hline  & &   & & & & & & & & 19 & 1 & 20  \\
\hline -1 & &   & & & & & & & & & 7 & -6  \\
\hline
\end{array}$
\end{center}
Now it is sufficient to look at the first row of $B$ to see that it is not an integer $^2\H(13;3)$.
\end{ex}

Now applying the previous corollary we obtain infinite classes of $\lambda$-fold relative Heffter arrays.
We start from the case $\lambda=2$ and $t=1$.

\begin{thm}\label{lambda2}
There exists a $^2\H(n;k)$ if and only if $n\geq k\geq 3$ with $nk\not\equiv 1 \pmod 4$.
\end{thm}
\begin{proof}
It is easy to see that a $^2\H(n;2)$ cannot exist since for every element $x$ of the array, the element $-x$
has to appear in the same row and in the same column of $x$, in order to satisfy condition ($\rm{c_1}$)
of Definition \ref{def:lambdaRelative}. So we can suppose $k\geq 3$.

By way of contradiction we suppose that there exists a $^2\H(n;k)$ with $nk\equiv 1 \pmod 4$.
This implies that $v=\frac{2nk}{\lambda}+t=nk+1\equiv 2 \pmod 4$. Since $\lambda=2$ and $t=1$ this cannot happen
in view of Proposition \ref{prop:nonexistence}.

Note that if $nk$ is even a $^2\H(n;k)$ exists in view of Theorem \ref{SMS} and Remark \ref{rem:magicHeffter}.
So we are left to consider the case $nk\equiv 3\pmod 4$.
It is easy to see an integer $\H(n;k)$ on $\Z_{2nk+1}$ can also be view as an integer $\H_2(n;k)$
on $\Z_{2nk+2}$. Thus by Theorem \ref{thm:squareinteger}, we have that there exists an
 integer $\H_2(n;k)$ if and only if $n\geq k\geq 3$ and $nk\equiv 0,3\pmod 4$.
 Hence the result follows by Corollary \ref{cor:proiezione}.
\end{proof}

For the case in which the array has an even number of filled cells in each column and
in each row, we present a complete solution both in the integer case and in the general one.

\begin{thm}\label{thm:skeveninteger}
Let $s,k$ be even positive integers.
There exists an integer\\ $^2\H(m,n;s,k)$ if and only if $ms=nk$, $2\leq s\leq n$, $2\leq k \leq m$
and
\begin{itemize}
\item[(1)] if $s=2$, either $n=2$ and $m=k\equiv 0\pmod 4$ or $n,k\geq3$;
\item[(2)] if $k=2$, either $m=2$ and $n=s\equiv 0\pmod 4$ or $m,s\geq 3$.
\end{itemize}
\end{thm}
\begin{proof}
Note that $ms=nk$, $2\leq s\leq n$, $2\leq k \leq m$ are the trivial necessary conditions for the existence of
a $^2\H(m,n;s,k)$. If $s=2$, by Remark \ref{rem:magicHeffter}, an integer $^2\H(m,n;2,k)$ is nothing but a $SMA(m,n;2,k)$
which exists if and only if either $n=2$ and $m=k\equiv 0,3\pmod 4$ or $n,k\geq3$ by Theorem \ref{thm:magic2}.
An analogous reasoning holds for $k=2$.

For $s$ and $k$ greater than $2$ the existence of an integer $^2\H(m,n;s,k)$ follows by Theorem \ref{thm:magicMP} and Remark \ref{rem:magicHeffter}.
\end{proof}

\begin{thm}\label{thm:skeven}
Let $s,k$ be even positive integers.
There exists a $^2\H(m,n;s,k)$ if and only if $ms=nk$, $2\leq s\leq n$, $2\leq k \leq m$ and $s,k$ are not both equal to $2$.
\end{thm}
\begin{proof}
We recall that $ms=nk$, $2\leq s\leq n$, $2\leq k \leq m$ are the trivial necessary conditions for the existence of
a $^2\H(m,n;s,k)$. Moreover, if $s=k=2$, then $n=m$ and we have seen in Theorem \ref{lambda2} that a $^2\H(n;2)$ does not exist.

On the other hand, in view of Theorem \ref{thm:skeveninteger},
to prove that the trivial necessary conditions are also sufficient, it remains
to show the existence of a
 $^2\H(m,2;2,k)$ with $m=k\equiv 2\pmod{4}$ and of a $^2\H(2,n;s,2)$
 with $n=s \equiv 2 \pmod4$. Note that in these cases the array has no empty cell.
Indeed,  because of the symmetry of Definition \ref{def:lambdaRelative},
 it suffices to prove the existence of a $^2\H(2,n;n,2)$ for any $n\equiv 2 \pmod{4}$ larger than $2$.
 At this purpose, for these values of $n$,  we consider the $2\times n$ array whose first row is:
$$\overline{R}_{1}=(-1,2,-3,4,-5,6,7,-8,-9,10,11,-12,\dots,-(4i+1),4i+2,4i+3,-(4i+4) ,$$
$$\ldots,-(n-5), n-4,n-3,-(n-2), (n-1),n)$$
and with $\overline{R}_2=-\overline{R}_1$.
We show that this array is a $^2\H(2,n;n,2)$.
Clearly  conditions $(\rm{a_1})$ and $(\rm{b_1})$ of Definition \ref{def:lambdaRelative} are satisfied.
Also, the sum of the elements in any column is $0$ because $\overline{R_2}=-\overline{R_1}$.
Let us consider the sum of the elements in the rows. Since $-(4i+1)+(4i+2)+(4i+3)-(4i+4)=0$, we have that:
$$\sum_{x\in\overline{R}_1} x=- \sum_{x\in\overline{R}_2} x=-1+2-3+4+(n-1)+n=2n+1.$$
This sum is zero modulo $v=2n+1$ and hence also the condition $(\rm{c_1})$ is satisfied.
\end{proof}

In the following we present a complete solution for the case in which the array has no empty cells,
such an array is also called \emph{tight}. Firstly we need to recall an elementary number theory lemma.
\begin{lem}\label{lem:requiredsum}
Let $s$ be an integer in the interval $\left[1,\frac{n(n+1)}{2}\right]$.
Then there exists a subset $S$ of $[1,n]$ whose elements sum to $s$, that is $\sum_{x\in S} x=s.$
\end{lem}
\begin{proof}
Let us suppose by contradiction there is no subset of $[1,n]$ that sums to $s$.
Let $S_1$ be a subset of $[1,n]$ whose sum $s_1$ is maximal among the ones smaller than $s$.
Let $y$ be the smallest element of $[1,n]$ that does not belong to $S_1$.
Clearly either $y=1$ or $y-1\in S_1$. We set $S_2$ to be $S_1\cup \{1\}$ in the first case and
$(S_1\setminus \{y-1\})\cup \{y\}$ in the second one. In both cases, $S_2$ sums to $s_1+1$ that is in contradiction with the maximality of $s_1$.
\end{proof}

\begin{thm} \label{thm:noemptycell}
There exists a $^2\H(m,n;n,m)$ if and only if $mn\not\equiv 1 \pmod 4$ and $m,n$ are not both equal to $2$.
\end{thm}
\begin{proof}
Firstly let us assume that $m=2$ (or, by symmetry, that $n=2$).
We note that a $^2\H(2;2)$ cannot exist. On the other hand, because of Theorem \ref{thm:skeven},
Theorem \ref{thm:SMS} and Remark \ref{rem:magicHeffter}, there exists a $^2\H(2,n;n,2)$ whenever
$n\not \equiv 1 \pmod{4}$ is larger than $2$. Thus, it suffices to construct a $^2\H(2,n;n,2)$ for any $n\equiv 1 \pmod{4}$.
Let us set $s=\frac{n^2-3n-2}{4}$. We note that, since $n\equiv 1 \pmod{4}$, $s$ is an integer.
We also have that $s$ is in the range $\left[1,\frac{n(n+1)}{2}\right]$. Therefore, due to Lemma \ref{lem:requiredsum},
there exists a subset $S$ of $[1,n]$ that sums to $s$.
Let now consider a $2\times n$ array $A$ such that:
\begin{itemize}
\item $\overline{R}_1$ contains the elements of $-S$ and the elements of $[1,n]\setminus S$;
\item $\overline{R}_2=-\overline{R}_1$.
\end{itemize}
Clearly $A$ satisfies conditions $(\rm{a_1})$ and $(\rm{b_1})$ of Definition \ref{def:lambdaRelative}.
Also the sum of the elements in any column is $0$ because $\overline{R_2}=-\overline{R_1}$.
Let us consider the sum of the elements in the rows:
$$\sum_{x\in\overline{R}_1} x=- \sum_{x\in\overline{R}_2} x=-\sum_{x\in S}x+\sum_{x\in [1,n]\setminus S}x=-2\sum_{x\in S}x+\sum_{x\in [1,n]}x=$$
$$-2\frac{n^2-3n-2}{4}+\frac{n^2+n}{2}=2n+1.$$
This sum is zero modulo $v=2n+1$ and hence also condition $(\rm{c_1})$ is satisfied.

Suppose now $m,n>2$.
A $^2\H(m,n;n,m)$ with $mn\equiv 1 \pmod 4$ does not exist in view of Proposition \ref{prop:nonexistence}.
The existence of a $^2\H(m,n;n,m)$ in the case $mn$ even follows by Theorem \ref{thm:SMS} and Remark \ref{rem:magicHeffter}.
When $mn\equiv 3 \pmod 4$, there exists an integer $\H(m,n;n,m)$ on $\Z_{2mn+1} $ as proved in Theorem 1.5 of
\cite{ABD}. This array can also be view as an integer $\H_2(m,n;n,m)$ on $\Z_{2mn+2}$.
So the existence of a $^2\H(m,n;n,m)$ follows by Corollary \ref{cor:proiezione}.
This completes the proof.
\end{proof}

Now we present some results that hold also for $\lambda>2$ and $t> 1$.

\begin{prop}\label{lambdak}
Let $n,k,\lambda$ be positive integers with $n\geq k\geq 3$ and $\lambda$ divisor of $k$.
  Then there exists a $^\lambda \H_{\frac{k}{\lambda}}(n;k)$ if one of the following is satisfied:
  \begin{itemize}
  \item[(1)] $k=5$ and $n\equiv 3 \pmod 4$;
  \item[(2)] $k\neq 5$ odd and $n\equiv 0,3 \pmod 4$.
  \end{itemize}
\end{prop}
\begin{proof}
The thesis follows by Theorem
\ref{thm:esistenzaRel} and Corollary \ref{cor:proiezione}.
\end{proof}

\begin{prop}\label{prop:n}
  For any odd integer $n\geq 3$ and for any divisor $\lambda$ of $n$ there exists a $^\lambda \H_\frac{n}{\lambda}(n;3)$.
\end{prop}
\begin{proof}
The result follows from Proposition \ref{n,2n} and Corollary \ref{cor:proiezione}.
\end{proof}

\begin{prop}\label{prop:2n}
  For any odd integer $n\geq 3$ and for any divisor $\lambda$ of $2n$ there exists a $^\lambda \H_\frac{2n}{\lambda}(n;3)$.
\end{prop}
\begin{proof}
The result follows from Proposition \ref{n,2n} and Corollary \ref{cor:proiezione}.
\end{proof}

\begin{prop}
Let $m,n,s,k,t,\lambda$ be positive integers such that $ms=nk$, $4\leq s\leq n$, $4\leq k\leq m$, $t$ is a divisor of $2nk$ and
$\lambda$ is a divisor of $t$. A $^\lambda \H_{\frac{t}{\lambda}}(m,n;s,k)$ exists in each of the following cases:
\begin{itemize}
  \item[(1)] $s\equiv k\equiv 0\pmod 4$;
  \item[(2)] either $s\equiv 2\pmod 4$ and $k\equiv 0\pmod 4$ or
 $s\equiv 0\pmod 4$ and $k\equiv 2\pmod 4$;
  \item[(3)] $s\equiv\ k \equiv 2\pmod 4$ and $m$ and $n$ are even.
\end{itemize}
\end{prop}

\begin{proof}
The result follows by Theorem \ref{main} and Corollary \ref{cor:proiezione}.
\end{proof}

Note that also by Theorems \ref{thm:magic2} and \ref{thm:magic3} one can get partial results about integer $2$-fold Heffter arrays,
other results about integer $\lambda$-fold relative Heffter arrays can be found in \cite{MP2}.

To conclude this section we present a general recursive construction.
\begin{prop}\label{prop:recursive}
If there exists an (integer) $\H_t(m,n;s,k)$ then there exists an (integer) $^{\alpha_1\lambda_1}\H_t(\lambda_1m,\lambda_2n;\alpha_1s,\alpha_2k)$
for all positive integers $\alpha_1\leq \lambda_2$, $\alpha_2\leq\lambda_1$ such that $\alpha_1\lambda_1=\alpha_2\lambda_2$,
$\alpha_1s\leq \lambda_2 n$ and $\alpha_2k\leq \lambda_1 m$.
\end{prop}
\begin{proof}
  Let $A=\H_t(m,n;s,k)$, so its elements belong to $\Z_{2nk+t}$. Set $B$ to be a $\lambda_1\times \lambda_2$ array with
  $\alpha_1$ filled cells in each row and $\alpha_2$ filled cells in each column.
  Clearly such an array exists since, by hypothesis, $\alpha_1\leq \lambda_2$, $\alpha_2\leq\lambda_1$ and $\alpha_1\lambda_1=\alpha_2\lambda_2$.
  It is easy to see that if we replace each filled cell of $B$ with the array $A$
  and each empty cell of $B$ with an empty $m\times n$ array we obtain
  a $^{\alpha_1\lambda_1}\H_t(\lambda_1m,\lambda_2n;\alpha_1s,\alpha_2k)$.
  Note that also the elements of the array so constructed belong to $\Z_{2nk+t}$,
  hence if $A$ is integer, then  the array $^{\alpha_1\lambda_1}\H_t(\lambda_1m,\lambda_2n;\alpha_1s,\alpha_2k)$
  is integer too.
\end{proof}

In view of the partial results presented in this section, we underline that the very general case in which the array is rectangular,
$\lambda$ and $t$ can assume any admissible value and at least one between $s$ and $k$ is odd is completely open.

\section{Relation with biembeddings}\label{sec:biembedding}

In \cite{A}, Archdeacon introduced Heffter arrays also because
 they are useful for finding biembeddings of cycle decompositions, as shown, for instance, in
\cite{CDDYbiem, CMPPHeffter, CPPBiembeddings, DM}.
In this section, generalizing some of his results we show how starting from a $\lambda$-fold relative Heffter array
it is possible to obtain suitable biembeddings.
In this context we equip a multigraph $\G$ with the following topology.
\begin{itemize}
\item If $\G$ is a simple graph, $\G$ is viewed with the usual topology as a $1$-dimensional simplicial complex.
\item If $\G$ is a multigraph, we consider the topology naturally induced on $\G$ by a simple graph $\G'$ that is a subdivision of $\G$.
\end{itemize}
Note that this topology is well defined because two multigraphs $\G$ and $\G'$ are homeomorphic if and only if there exists an isomorphism from some subdivision of $\G$ to some subdivision of $\G'$. Now we provide the following definition, see \cite{Moh} for the simple graph case.
\begin{defi}
An \emph{embedding} of a multigraph $\G$ in a surface $\Sigma$ is a continuous injective
mapping $\psi: \G \to \Sigma$, where $\G$ is viewed with the topology described above.
\end{defi}

The  connected components of $\Sigma \setminus \psi(\G)$ are called $\psi$-\emph{faces}.
If each $\psi$-face is homeomorphic to an open disc, then the embedding $\psi$ is said to be \emph{cellular}.

\begin{defi}
 A \emph{biembedding} of  two cycle decompositions $\D$ and $\D'$ of a multigraph $\G$  is a face $2$-colorable embedding
  of $\G$ in which one color class is comprised of the cycles in $\D$ and  the other class contains the cycles in $\D'$.
\end{defi}

If we are not interested in the decompositions, we will simply speak of biembeddings of multigraphs.

Following the notation given in \cite{A}, for every edge $e$ of a  multigraph $\G$, let $e^+$ and $e^-$ denote
its two possible directions
and let $\tau$ be the involution swapping $e^+$ and $e^-$ for every $e$.
Let $D(\G)$ be the set of all directed edges of $\G$ and, for any $v\in V(\G)$, call $D_v$ the set of edges directed out
of $v$.
A local rotation $\rho_v$ is a cyclic permutation of $D_v$. If we select a local rotation for each vertex of $\G$, then
all together
they form a rotation of $D(\G)$. We recall the following result, see \cite{A, GT, MT}.

\begin{thm}\label{thm:embedding}
 A rotation $\rho$ on a multigraph $\G$ is equivalent to a cellular embedding of $\G$ in an orientable surface.
 The face boundaries of the embedding corresponding to $\rho$ are the orbits of $\rho \circ \tau$.
\end{thm}

Given a  $\lambda$-fold relative Heffter array $^\lambda\H_t(m,n; s,k)$, say $A$, then the orderings $\omega_r$ (associated to the permutation $\alpha_r$) and $\omega_c$ (associated to the permutation $\alpha_c$) are said to be
\emph{compatible} if $\alpha_c \circ \alpha_r$ is a cycle of length $|skel(A)|$.
With the same proof of Theorem 1.4 of \cite{CDDYbiem} we obtain the following necessary conditions for the existence of compatible orderings.
\begin{prop}
If there exist compatible orderings $\omega_r$ and $\omega_c$ for a $\lambda$-fold relative Heffter array $^\lambda\H_t(m,n; s,k)$ then one of the following has to be satisfied:
\begin{itemize}
\item[(1)] $m,n,s$ and $k$ are odd;
\item[(2)] $m$ is odd, while $n$ and $s$ are even;
\item[(3)] $n$ is odd, while $m$ and $k$ are even.
\end{itemize}
\end{prop}
\begin{thm}\label{thm:biembedding}
Let $A$ be a  $\lambda$-fold relative Heffter array $^\lambda\H_t(m,n; s,k)$ that is simple with respect to the compatible orderings $\omega_r$ and $\omega_c$.
Then there exists a cellular biembedding of the cyclic cycle decompositions $\mathcal{D}_{\omega_r^{-1}}$ and
$\mathcal{D}_{\omega_c}$ of $^\lambda K_{(\frac{2nk}{\lambda t}+1)\times t}$ into an orientable surface.
\end{thm}

\begin{proof}
Since the orderings $\omega_r$ and $\omega_c$ are compatible, we have that $\alpha_c\circ\alpha_r$ is a cycle of length
$|skel(A)|$. Let us consider the permutation $\gamma$ on the set $skel(A)\times \{1,-1\}$ defined by:
$$\gamma(i,j,\varepsilon)=\begin{cases}
(\alpha_r(i,j),-1)\mbox{ if } \varepsilon=1;\\
(\alpha_c(i,j),1)\mbox{ if } \varepsilon=-1.\\
\end{cases}$$
Note that, if $(i,j)\in skel(A)$, then $\gamma^2(i,j,1)=(\alpha_c\circ\alpha_r(i,j),1)$ and hence $\gamma^2$ acts
cyclically on $skel(A)\times \{1\}$. Also $\gamma$ exchanges $skel(A)\times \{1\}$ with $skel(A)\times \{-1\}$. Thus it acts cyclically on $skel(A)\times \{1,-1\}$. We also define $\gamma_1$ to be the projection of $\gamma$ on $skel(A)$ and by $\gamma_2$ to be its projection on $\{1,-1\}$.

We note that the oriented edges of the multigraph $^\lambda K_{(\frac{2nk}{\lambda t}+1)\times t}$
having $x$ in the first component are of the form
$(x,x\pm a_{i,j})$ where we consider the edges $(x,x+\varepsilon\cdot a_{i,j})$ and $(x,x+\varepsilon '\cdot a_{i',j'})$
to be equal only if $(i,j,\varepsilon)=(i',j',\varepsilon ')$.
In particular if  $a_{i,j}=a_{i',j'}$ but $(i,j)\not=(i',j')$, then $(x,x+\varepsilon \cdot a_{i,j})$ and $(x,x+\varepsilon \cdot a_{i',j'})$ are different edges of the multigraph $^\lambda K_{(\frac{2nk}{\lambda t}+1)\times t}$.
Therefore we define the map $\rho_x$ on those edges so that:
$$\rho_x((x,x+\varepsilon\cdot a_{i,j}))= (x,x+\gamma_2(i,j,\varepsilon)\cdot a_{\gamma_1(i,j,\varepsilon)}).$$
Since $\gamma$ acts cyclically on $skel(A)\times \{1,-1\}$ the map $\rho_x$ is a local rotation at $x$ and we can define a rotation $\rho$ on the edges of the multigraph $^\lambda K_{(\frac{2nk}{\lambda t}+1)\times t}$ so that $\rho((x,x+a_{i,j}))=\rho_x((x,x+a_{i,j}))$.
Hence, by Theorem \ref{thm:embedding}, there exists a cellular embedding $\sigma$ of  $^\lambda K_{(\frac{2nk}{\lambda t}+1)\times t}$ in an orientable surface so that the face boundaries correspond to the orbits of $\rho\circ \tau$ where $\tau((x,x\pm a_{i,j}))=(x\pm a_{i,j},x)$.
Let us consider the oriented edge $(x,x+a_{i,j})$ with $a_{i,j} \in \E(A)$, and let  $\overline{C}$ be the column containing $a_{i,j}$.
We have that:
$$\rho\circ\tau((x,x+a_{i,j}))=\rho((x+a_{i,j},(x+a_{i,j})-a_{i,j}))=(x+a_{i,j},(x+a_{i,j})+a_{\alpha_c(i,j)}).$$
Thus $(x,x+a_{i,j})$ belongs to the boundary of the face $F_1$ delimited by the oriented edges:
$$(x,x+a_{i,j}),(x+a_{i,j},x+a_{i,j}+a_{\alpha_c(i,j)}),(x+a_{i,j}+a_{\alpha_c(i,j)},x+a_{i,j}+a_{\alpha_c(i,j)}+a_{\alpha_c^2(i,j)}),\dots$$
$$\dots ,\left(x+\sum_{t=0}^{|\E(\overline{C})|-2} a_{\alpha_c^t(i,j)},x\right).$$
We note that the cycle associated to the face $F_1$ is:
$$\left(x,x+a_{i,j},x+a_{i,j}+a_{\alpha_c(i,j)},\ldots,x+\sum_{t=0}^{|\E(\overline{C})|-2} a_{\alpha_c^t(i,j)}\right).$$
Let us now consider the oriented edge $(x,x-a_{i,j})$ with $a_{i,j} \in \E(A)$ and let us name by $\overline{R}$ the row containing the element $a_{i,j}$. We have
that:
$$\rho\circ\tau((x,x-a_{i,j}))=\rho((x-a_{i,j},(x-a_{i,j})+a_{i,j}))=(x-a_{i,j},(x-a_{i,j})-a_{\alpha_r(i,j)}).$$
Thus $(x,x-a_{i,j})$ belongs to the boundary of the face $F_2$ delimited by the oriented edges:
$$(x,x-a_{i,j}),(x-a_{i,j},x-a_{i,j}-a_{\alpha_r(i,j)}),$$
$$(x-a_{i,j}-a_{\alpha_r(i,j)},x-a_{i,j}-a_{\alpha_r(i,j)}-a_{\alpha^2_r(i,j)}),\dots, \left(x-\sum_{t=0}^{|\E(\overline{R})|-2}
a_{\alpha^t_r(i,j)},x\right).$$
Since $A$ is a $\lambda$-fold relative Heffter array and $\alpha_r$ acts cyclically on $\E(\overline{R})$, for any $j\in [1,
|\E(\overline{R})|]$ we have that:
$$-\sum_{t=0}^{j-1} a_{\alpha^t_r(i,j)}=
\sum_{t=j}^{|\E(\overline{R})|-1}a_{\alpha_r^t(i,j)}=\sum_{t=1}^{|\E(\overline{R})|-j}a_{\alpha_r^{|\E(\overline{R})|-t}(i,j)}=\sum_{t=1}^{|\E(\overline{R})|-j}a_{\alpha_r^{-t}(i,j)}.$$
It follows that the cycle associated to the face $F_2$ can be written also as:
$$\left(x,x+\sum_{t=1}^{|\E(\overline{R})|-1}a_{\alpha_r^{-t}(i,j)},x+\sum_{t=1}^{|\E(\overline{R})|-2}a_{\alpha_r^{-t}(i,j)},
\dots,x+a_{\alpha_r^{-1}(i,j)}\right).$$
Therefore any nonoriented edge $\{x,x-a_{i,j}\}$ belongs to the  boundaries of exactly two faces: one of type $F_1$ and one of
type $F_2$. Hence the embedding is 2-colorable.

Moreover, it is easy to see that those face boundaries are the cycles obtained from the $\lambda$-fold relative Heffter array $A$ following the
orderings  $\omega_c$ and $\omega_r^{-1}$.
\end{proof}

As already remarked in \cite{CPPBiembeddings},
looking for compatible orderings in the case of a globally simple $\lambda$-fold Heffter array led us to investigate the following problem
introduced in \cite{CDP}.
Let $A$ be an $m\times n$ \emph{toroidal} p.f. array. By $r_i$ we denote the orientation of the $i$-th row,
precisely $r_i=1$ if it is from left to right and $r_i=-1$ if it is from right to left. Analogously, for the $j$-th
column, if its orientation $c_j$ is  from  top to bottom then $c_j=1$ otherwise $c_j=-1$. Assume that an orientation
$\R=(r_1,\dots,r_m)$
 and $\C=(c_1,\dots,c_n)$ is fixed. Given an initial filled cell $(i_1,j_1)$ consider the sequence
$ L_{\R,\C}(i_1,j_1)=((i_1,j_1),(i_2,j_2),\ldots,(i_\ell,j_\ell),$ $(i_{\ell+1},j_{\ell+1}),\ldots)$
where $j_{\ell+1}$ is the column index of the filled cell $(i_\ell,j_{\ell+1})$ of the row $\overline{R}_{i_\ell}$ next
to
$(i_\ell,j_\ell)$ in the orientation $r_{i_\ell}$,
and where $i_{\ell+1}$ is the row index of the filled cell of the column $\overline{C}_{j_{\ell+1}}$ next to
$(i_\ell,j_{\ell+1})$ in the orientation $c_{j_{\ell+1}}$.
The problem is the following:

\begin{KN}
Given a toroidal p.f. array $A$,
do there exist $\R$ and $\C$ such that the list $L_{\R,\C}$ covers all the filled
cells of $A$?
\end{KN}

By $P(A)$ we will denote the \probname\ for a given array $A$.
Also, given a filled cell $(i,j)$, if $L_{\R,\C}(i,j)$ covers all the filled positions of $A$ we will
say that  $(\R,\C)$ is a solution of $P(A)$.
For known results about this problem see \cite{CDP}.
The relationship between the Crazy Knight's Tour Problem and globally simple $\lambda$-fold Heffter arrays is explained in the following result
which is an easy consequence of Theorem
\ref{thm:biembedding}.

\begin{cor}\label{preprecedente}
  Let $A$ be a globally simple $^{\lambda}\H_t(m,n;s,k)$  such that $P(A)$ admits a solution $(\R,\C)$.
  Then there exists an orientable biembedding of  $^\lambda K_{(\frac{2nk}{\lambda t}+1)\times t}$ such that every edge is on a face
  of size $s$ and a face of size $k$.
\end{cor}
Therefore the first ingredient we need in order to construct biembeddings is a globally simple $\lambda$-fold
relative Heffter array.
Clearly, if $k=3$, then any ordering is simple and hence also the natural one.
On the other hand for larger values of $k$ the simplicity condition is not necessarily
satisfied by the natural orderings.  It is easy to see that when either $k=4$ or $k=5$ a $^\lambda\H_t(n;k)$ is globally simple if and only if it
does not contain two opposite elements consecutive in a row or in a column.
For example, one can check that, when $n\equiv 3 \pmod{4}$, the $^{2}\H(n;5)$
constructed via Theorem \ref{lambda2} is not globally simple because in the last row we have two consecutive
cells filled with $-2n$ and $2n$.
Although in general it seems quite difficult to obtain a $^\lambda \H_t(n;k)$ with this property, here we are able to provide a direct
construction of a globally simple $^{2}\H(n;5)$ whenever $n\equiv 3\pmod{4}$.
In order to present this construction we introduce some notation.
Let $A$ be an $n\times n$ p. f. array, for $i\in [1,n]$ the $i$-th diagonal is so defined
$D_i=\{(i,1),(i+1,2),\ldots,(i-1,n)\}$. All the arithmetic on the row and column indices is
performed modulo $n$, where the set of reduced residues is $\{1,2,\ldots,n\}$.
We say that the diagonals $D_i, D_{i+1},\ldots, D_{i+k}$ are $k+1$ consecutive diagonals.
Given a positive integer $k$ we say that a square Heffter array $A$ of size $n\geq k$ is \emph{cyclically $k$-diagonal} if
the nonempty cells of $A$ are exactly those of $k$ consecutive diagonals.
For instance the array of Example \ref{ex:inv} is cyclically $3$-diagonal.

\begin{rem}\label{rem:ciclici}
We point out that, if $k=3$, the arrays of Theorems \ref{thm:squareinteger} and \ref{thm:esistenzaRel} are
cyclically $3$-diagonal. This implies that if $k=3$ and $n$ is odd, also the arrays of Theorem \ref{lambda2}
and Proposition \ref{lambdak}
are cyclically $3$-diagonal.

The same holds for the arrays of Proposition \ref{n,2n} and hence also those of Propositions \ref{prop:n}
and \ref{prop:2n}
have the same property.
\end{rem}

Now we recall the procedure introduced in \cite{DW}, since it is very useful to describe our construction.
In an $n\times n$ array $A$ the procedure $diag(r,c,s,\Delta_1,\Delta_2,\ell)$ installs the entries
$$A[r+i\Delta_1,c+i\Delta_1]=s+i\Delta_2\ \textrm{for}\ i\in[0,\ell-1],$$
where by $A[i,j]$ we mean the element of $A$ in position $(i,j)$.
The parameters used in the $diag$ procedure have the following meaning:
\begin{itemize}
  \item $r$ denotes the starting row,
  \item $c$ denotes the starting column,
  \item $s$ denotes the entry $A[r,c]$,
  \item $\Delta_1$ denotes the increasing value of the row and column at each step,
  \item $\Delta_2$ denotes how much the entry is changed at each step,
  \item $\ell$ is the length of the chain.
\end{itemize}

\begin{prop}\label{prop:5GS}
Let $n\equiv 3 \pmod 4$ with $n\geq7$. Then  there exists a cyclically $5$-diagonal globally simple $^2\H(n;5)$.
\end{prop}
\begin{proof}
Let $H=(h_{i,j})$ be the integer $\H(n;5)$ described in \cite{ADDY} for $n\equiv 3\pmod{4}$ and define $a_{2i,2j}=h_{i,j}\pmod{5n+1}$.
Then $A=(a_{i,j})$ is a cyclically $5$-diagonal $^2\H(n;5)$ since the non-empty
cells of $H$ are exactly those of the diagonals $D_1$, $D_2$, $D_{\frac{n+1}{2}}, D_{\frac{n+3}{2}}, D_n$.
One can check, by a simple but quite long computation, that $A$ can be written also using the following procedures labeled $\texttt{A}$ to  $\texttt{N}$ where the sums are considered modulo $5n+1$.
$$\begin{array}{rlcrl}
\texttt{A}: &  diag\left(3,3,\frac{n-3}{2},2,-1,\frac{n-5}{2}\right); & \hfill &
\texttt{B}: &  diag\left(4,4,-(n-2),2,1,\frac{n-3}{2}\right);\\[3pt]
\texttt{C}: &  diag\left(3,2,n+1,2,2,\frac{n-1}{2}\right); &&
\texttt{D}: &  diag\left(4,3,3n-1,2,-2,\frac{n-3}{2}\right);\\[3pt]
\texttt{E}: &  diag\left(2,3,-3n,2,2,\frac{n-1}{2}\right); &  &
\texttt{F}: &  diag\left(3,4,-(n+2),2,-2,\frac{n-3}{2}\right);\\[3pt]
\texttt{G}: &  diag\left(3,1,-\frac{15n+3}{4},4,1,\frac{n-3}{4}\right); &&
\texttt{H}: &  diag\left(4,2,-(3n+3),4,-1,\frac{n+1}{4}\right);\\[3pt]
\texttt{I}: &  diag\left(5,3,-\frac{19n-9}{4},4,1,\frac{n-3}{4}\right);&&
\texttt{J}: &  diag\left(6,4,-(4n+1),4,-1,\frac{n-3}{4}\right);\\[3pt]
\texttt{K}: &  diag\left(1,3,\frac{17n+1}{4},4,1,\frac{n-3}{4}\right); &&
\texttt{L}: &  diag\left(2,4,5n-2,4,-1,\frac{n+1}{4}\right);\\[3pt]
\texttt{M}: &  diag\left(3,5,\frac{13n+13}{4},4,1,\frac{n-3}{4}\right); &&
\texttt{N}: &  diag\left(4,6,4n,4,-1,\frac{n-3}{4}\right).
\end{array}$$
We also fill the following cells in an \textit{ad hoc} manner:
$$\begin{array}{lclcl}
A\left[1,1\right]=-n, & \hfill &  A\left[1,2\right]=-2n+1, &  \hfill & A\left[1,n\right]=2n+1,\\[3pt]
A\left[2,1\right]=2n+2, & & A\left[2,2\right]=n-1, &  & A\left[2,n\right]=-(5n-1),\\[3pt]
A\left[n-2,n-2\right]=-\frac{n-1}{2}, & & A\left[n-2,n\right]=5n, && A\left[n,1\right]=-2n,\\[3pt]
A\left[n,2\right]=3n+2, & & A\left[n,n-2\right]=-(3n+1), & & A\left[n,n\right]=1.
\end{array}$$

We now prove that the array $A$ is globally simple.
To aid in the proof we give a schematic picture of where each of the diagonal procedures fills cells (see Figure \ref{fig4}).
We have placed an $\texttt{X}$ in the \textit{ad hoc} cells.
We list the elements in every row.

\begin{figure}
 \begin{footnotesize}
\begin{center}
\begin{tabular}{|c|c|c|c|c|c|c|c|c|c|c|c|c|c|c|}
  \hline
 \texttt{X} & \texttt{X} & \texttt{K} & & & & & & & & & & & \texttt{H}&  \texttt{X}\\ \hline
 \texttt{X} & \texttt{X} & \texttt{E} &\texttt{L} & & & & & & & & & & &\texttt{X}  \\ \hline
\texttt{G} &\texttt{C} & \texttt{A} & \texttt{F} &\texttt{M} & & & & & & & & & &   \\ \hline
 & \texttt{H}&\texttt{D} & \texttt{B} & \texttt{E}  &\texttt{N} & & & & & & & & &   \\ \hline
 & &\texttt{I} &\texttt{C} & \texttt{A} & \texttt{F} &\texttt{K} & & & & & & & &    \\ \hline
& & & \texttt{J}&\texttt{D} & \texttt{B} & \texttt{E}  &\texttt{L} & & & & & & &   \\ \hline
& & & &\texttt{G} &\texttt{C} & \texttt{A} & \texttt{F} &\texttt{M} & & & & & &   \\ \hline
& & & & &\texttt{H} &\texttt{D} & \texttt{B} & \texttt{E} &\texttt{N} & & & & &   \\ \hline
& & & & & & \texttt{I}&\texttt{C} & \texttt{A} & \texttt{F} &\texttt{K} & & & &    \\ \hline
& & & & & & & \texttt{J}&\texttt{D} & \texttt{B} & \texttt{E} &\texttt{L} & & &    \\ \hline
& & & & & & & &\texttt{G} &\texttt{C} & \texttt{A} & \texttt{F} &\texttt{M} & &     \\ \hline
& & & & & & & & & \texttt{H}&\texttt{D} & \texttt{B} & \texttt{E} &\texttt{N} &     \\ \hline
& & & & & & & & & &\texttt{I} &\texttt{C} & \texttt{X} & \texttt{F} & \texttt{X}    \\ \hline
\texttt{L}& & & & & & & & & & & \texttt{J}&\texttt{D} & \texttt{B} & \texttt{E}    \\ \hline
\texttt{X} &\texttt{X} & & & & & & & & & & &\texttt{X} & \texttt{C} & \texttt{X}     \\ \hline
\end{tabular}
\end{center}
\end{footnotesize}
\caption{Scheme of construction with $n=15$.}\label{fig4}
\end{figure}

\begin{description}
\item[Row $1$] There are three \textit{ad hoc} values, the first  of the \texttt{K} diagonal and the last  of the
\texttt{H} diagonal. Namely:  $$\E(\overline{R}_1)=(-n,-2n+1,\frac{17n+1}{4},-\frac{13n+9}{4},2n+1).$$
\item[Row $2$] There are three \textit{ad hoc} values, the first of the \texttt{E} diagonal as well as the first of
the \texttt{L} diagonal. Hence $$\E(\overline{R}_2)=(2n+2,n-1,-3n,5n-2,-(5n-1)).$$
\item[Row $3$ to $n-3$] Consider the row $r$, there are four cases according to the congruence class of $r$ modulo $4$.
If $r\equiv 3 \pmod 4$ write $r=3+4i$ where $i\in \left[0, \frac{n-7}{4}\right]$. Notice that from the \texttt{G},
\texttt{C}, \texttt{A}, \texttt{F} and \texttt{M} diagonal cells we get the following list:
$$\E(\overline{R}_{3+4i})=\left (-\frac{15n+3}{4}+i,n+1 +4i,\frac{n-3}{2}-2i,-(n+2+4i),\frac{13n+13}{4}+i\right).$$
If $r\equiv 0 \pmod 4$ write $r=4+4i$ where $i\in \left[0, \frac{n-7}{4}\right]$. Notice that from the \texttt{H},
\texttt{D}, \texttt{B}, \texttt{E} and \texttt{N} diagonal cells we get the following list:
$$\E(\overline{R}_{4+4i})=\left(-(3n+3+i),3n-1 -4i,-(n-2-2i),-(3n-2-4i),4n-i\right).$$
If $r\equiv 1 \pmod 4$ write $r=5+4i$ where $i\in \left[0, \frac{n-11}{4}\right]$. Notice that from the \texttt{I},
\texttt{C}, \texttt{A}, \texttt{F} and \texttt{K} diagonal cells we get the following list:
$$\E(\overline{R}_{5+4i})=\left(-\frac{19n-9}{4}+i,n+3+4i,\frac{n-5}{2}-2i,-(n+4+4i),\frac{17n+5}{4}+i\right).$$
If $r\equiv 2 \pmod 4$ write $r=6+4i$ where $i\in \left[0, \frac{n-11}{4}\right]$. Notice that from the \texttt{J},
\texttt{D}, \texttt{B}, \texttt{E} and \texttt{L} diagonal cells we get the following list:
$$\E(\overline{R}_{6+4i})=\left(-(4n+1+i),3n-3 -4i,-(n-3-2i),-(3n-4-4i),5n-3-i\right).$$
\item[Row $n-2$] This row contains two \textit{ad hoc} values, the last of the \texttt{I} diagonal, the
$\frac{n-3}{2}$-th element of the \texttt{C} diagonal and the last of \texttt{F} diagonal. Hence
    $$\E(\overline{R}_{n-2})=\left(-\frac{9n-1}{2},2n-4,-\frac{n-1}{2},-(2n-3),5n\right).$$
\item[Row $n-1$] We have the last elements of  the \texttt{L}, \texttt{J},  \texttt{D},  \texttt{B} and \texttt{E}
diagonals. Namely $$\E(\overline{R}_{n-1})=\left(\frac{19n-5}{4},-\frac{17n-3}{4},
2n+4,-\frac{n+1}{2},-(2n+3)\right).$$
\item[Row $n$] This row contains four \textit{ad hoc} values and the last of \texttt{C} diagonal. The list is
$$\E(\overline{R}_{n})=\left(-2n,3n+2,-(3n+1),2n-2,1\right).$$
   \end{description}
It is easy to see that in each row we do not have any consecutive pair of type $x,$ $-x$ (modulo $5n+1$),
hence the rows are simple with respect to the natural ordering.
 Now we check that also the columns have the same property.
\begin{description}
\item[Column $1$] There are three \textit{ad hoc} values, the first  of the \texttt{G} diagonal and the last  of the
\texttt{L} diagonal. Namely
$$\E(\overline{C}_{1})=\left(-n,2n+2,-\frac{15n+3}{4},\frac{19n-5}{4},-2n\right).$$
\item[Column $2$] There are three \textit{ad hoc} values, the first of the \texttt{C} diagonal as well as the first
of the \texttt{H} diagonal. The list is $$\E(\overline{C}_{2})=\left(-2n+1,n-1,n+1,-(3n+3),3n+2\right).$$
\item[Column $3$ to $n-3$] Consider the column $c$, there are four cases according to the congruence class of $c$ modulo
$4$.
If $c\equiv 3 \pmod 4$ write $c=3+4i$ where $i\in \left[0, \frac{n-7}{4}\right]$. Notice that from the \texttt{K},
\texttt{E}, \texttt{A}, \texttt{D} and \texttt{I} diagonal cells we get the following list:
$$\E(\overline{C}_{3+4i})=\left(\frac{17n+1}{4}+i,-(3n-4i),\frac{n-3}{2}-2i,3n-1-4i,-\frac{19n-9}{4}+i\right).$$
If $c\equiv 0 \pmod 4$ write $c=4+4i$ where $i\in \left[0, \frac{n-7}{4}\right]$. Notice that from the \texttt{L},
\texttt{F}, \texttt{B}, \texttt{C} and \texttt{J} diagonal cells we get the following list:
$$\E(\overline{C}_{4+4i})=\left(5n-2-i,-(n+2+4i),-(n-2-2i),n+3+4i,-(4n+1+i)\right).$$
If $c\equiv 1 \pmod 4$ write $c=5+4i$ where $i\in \left[0, \frac{n-11}{4}\right]$. Notice that from the \texttt{M},
\texttt{E}, \texttt{A}, \texttt{D} and \texttt{G} diagonal cells we get the following list:
$$\E(\overline{C}_{5+4i})=\left(\frac{13n+13}{4}+i,-(3n-2-4i),\frac{n-5}{2}-2i,3n-3-4i,-\frac{15n-1}{4}+i\right).$$
If $c\equiv 2 \pmod 4$ write $c=6+4i$ where $i\in \left[0, \frac{n-11}{4}\right]$. Notice that from the \texttt{N},
\texttt{F}, \texttt{B}, \texttt{C} and \texttt{H} diagonal cells we get the following list:
$$\E(\overline{C}_{6+4i})=\left(4n-i,-(n+4+4i),-(n-3-2i),n+5+4i,-(3n+4+i)\right).$$
\item[Column $n-2$] This column contains two \textit{ad hoc} values, the last of the \texttt{M} diagonal, the
$\frac{n-3}{2}$-th element of the \texttt{E} diagonal and the last of \texttt{D} diagonal. The list is
    $$\E(\overline{C}_{n-2})=\left(\frac{7n+3}{2},-(2n+5),-\frac{n-1}{2},2n+4,-(3n+1)\right).$$
\item[Column $n-1$] There are the last elements of  the \texttt{H}, \texttt{N},  \texttt{F},  \texttt{B} and \texttt{C}
diagonals. Hence
$$\E(\overline{C}_{n-1})=\left(-\frac{13n+9}{4},\frac{15n+7}{4},-(2n-3),-\frac{n+1}{2},2n-2\right).$$
\item[Column $n$] This column contains four \textit{ad hoc} values and the last of \texttt{E} diagonal. The list is
$$\E(\overline{C}_{n})=\left(2n+1,-(5n-1),5n,-(2n+3),1\right).$$
\end{description}
It is easy to see that in each column we do not have any consecutive pair of type $x,$  $-x$ (modulo $5n+1$),
hence also the columns are simple with respect to the natural ordering.
Thus, $A$ is a globally simple cyclically $5$-diagonal $^2\H(n;5)$ for every $n\equiv 3 \pmod 4$.
\end{proof}
\begin{ex}
Following the proof of Proposition \ref{prop:5GS} we obtain the
globally simple $^2\H(15;5)$ below.
\begin{center}
\begin{tiny}
$\begin{array}{|r|r|r|r|r|r|r|r|r|r|r|r|r|r|r|}
\hline -15 & -29 & 64  &  &  &  &  & & & &&&& -51 & 31  \\
\hline  32 & 14 & -45 & 73 &  &  &  &  & & &&&&& -74 \\
\hline  -57 & 16 & 6 & -17 & 52  &  &  & & &&&&& & \\
\hline   & -48  & 44 & -13 & -43 & 60 &  & &&&& & & &\\
\hline   &  & -69 & 18   & 5 & -19  & 65 & & & & & & & &\\
\hline   &  &  & -61  &  42 & -12 & -41 & 72 & & & & & & & \\
\hline    &  &  &  & -56 & 20 & 4 & -21 & 53 & & & & & & \\
\hline   &  &  &  &  & -49  & 40 & -11  & -39 & 59& & & & &\\
\hline   &  &  &   &  &  & -68 & 22 & 3 & -23 & 66 & & & & \\
\hline   &  &  &   &  &  & &-62 & 38 & -10 & -37 & 71  & & & \\
\hline   &  &  &   &  &  & & &-55 & 24 & 2 & -25 & 54 & &  \\
\hline   &  &  &   &  &  & & & &-50 & 36 & -9 & -35 & 58  & \\
\hline   &  &  &   &  &  & & & & &-67 & 26 & -7 & -27 & 75 \\
\hline   70  &  &   &  &  & & & & & & &-63 & 34 & -8 & -33 \\
\hline   -30& 47  &  &   &  &  & & & & & &  & -46 & 28 & 1 \\
\hline
\end{array}$
\end{tiny}
\end{center}
\end{ex}

Now we are able to give infinite classes of biembeddings.
\begin{prop}
There exists a cellular biembedding of a pair of cyclic $k$-cycle decompositions of $^{\lambda}K_{\frac{2nk}{\lambda}+1}$ into an orientable surface in each of the following cases:
\begin{itemize}
\item[(1)] $k=3$, $\lambda=2$ and $n\equiv 1 \pmod{4}$;
\item[(2)] $k=3$, $\lambda=3$ and $n\equiv 3 \pmod{4}$;
\item[(3)] $k=5$, $\lambda=2$ and $n\equiv 3 \pmod{4}$.
\end{itemize}
\end{prop}
\begin{proof}
The result follows from Theorem \ref{thm:biembedding}
once we have provided a globally simple $\lambda$-fold Heffter array with compatibile natural orderings.
\begin{itemize}
\item[(1,2)] For those values of $n$, a $^2\H(n;3)$ and a $^3\H(n;3)$ are constructed respectively in Theorem \ref{lambda2} and in
Proposition \ref{lambdak}. Since they are cyclically $3$-diagonal by Remark \ref{rem:ciclici} the compability
of their natural orderings follows from  Proposition 3.4 of \cite{CMPPHeffter}.
\item[(3)] Let $n\equiv 3 \pmod{4}$, a globally simple $^2\H(n;5)$ is given in Proposition \ref{prop:5GS}.
Since these arrays are cyclically $5$-diagonal, as before, the compability of their natural orderings
     follows from  Proposition 3.4 of \cite{CMPPHeffter}.
\end{itemize}
\end{proof}

\begin{prop}
  There exists a cellular biembedding of a pair of cyclic $3$-cycle decompositions of $^{\lambda}K_{(\frac{6n}{ t}+1)\times \frac{t}{\lambda}}$ into an orientable surface when $n$ is odd, $t=n,2n$ and $\lambda$ divides $t$.
\end{prop}
\begin{proof}
Reasoning as in the proof of previous proposition, one can see that
the result follows from Theorem \ref{thm:biembedding}, Propositions \ref{prop:n} and \ref{prop:2n}, Remark \ref{rem:ciclici}
and Proposition 3.4 of \cite{CMPPHeffter}.
\end{proof}

\section{Some topological considerations}
In this paper we have introduced the concept of relative Heffter arrays and
we have provided constructions for infinite families of such objects.
The main tool we have used is Theorem \ref{thm:proiezione} that allows us to
obtain an $^{\alpha\lambda} \H_{\frac{t}{\lambda}}(m,n;s,k)$, say $B$, from an $^{\alpha}\H_{t}(m,n;s,k)$, say $A$.
Then, in the previous section, we have seen that if $B$ admits simple and compatible orderings
$\omega_r$ and $\omega_c$, there exists a cellular biembedding $\sigma$ of
$^{\alpha\lambda} K_{(\frac{2nk}{\alpha t}+1)\times \frac{t}{\lambda}}$ into an orientable surface $\Sigma$.
If $\omega_r$ and $\omega_c$ are simple and compatible the orderings $\omega_r'$ and $\omega_c'$ of $A$
associated to the same permutations $\alpha_r$ and $\alpha_c$ are simple and compatible too:
in fact the compatibility only depends on the skeleton of $A$ that is the same of $B$
and the simplicity follows from that of $\omega_r$ and $\omega_c$. This means that
there exists also a cellular biembedding $\sigma'$ of $^\alpha K_{(\frac{2nk}{\alpha t}+1)\times t}$
into an orientable surface $\Sigma'$. In this section we want to study the relations between $\sigma$
and $\sigma'$ and between $\Sigma$ and $\Sigma'$. At this purpose we recall the following definition.

\begin{defi}
Let $\Sigma$ be a topological space. A \emph{covering space} of $\Sigma$ is a topological space $\Sigma'$
together with a continuous surjective map $p:\Sigma'\rightarrow \Sigma$ with the following property: for every $x\in \Sigma$,
there exists an open neighborhood $U$ of $x$ such that $p^{-1}(U)$ is a union of disjoint open sets
in $\Sigma'$ each of which is mapped homeomorphically onto $U$ by $p$.
The map $p$ is called \emph{covering map}.
\end{defi}
If we consider a (multi)graph $\Gamma$ as a topological space, a (multi)graph $\Gamma'$ is a covering space for $\Gamma$
if there exists a surjective map $\pi: \Gamma'\rightarrow \Gamma$ such that maps edges incident to $x\in V(\Gamma')$ one-to-one onto edges incident to $\pi(x)$. In this case $\Gamma'$ is also called \emph{covering (multi)graph}.

\begin{ex}\label{ex:projection}
Let us consider the multipartite graph $K_{\frac{v'}{t}\times t}$. We can identify its vertex set with $\mathbb{Z}_{v'}$ and its edges with the pairs $\{x,x+a\}$ such that $x\in \mathbb{Z}_{v'}$ and $a$ is not in $\frac{v'}{t}\mathbb{Z}_{v'}$.
Now, given a divisor $\lambda$ of $t$ we can consider the multipartite multigraph
$^{\lambda} K_{\frac{v'}{t}\times \frac{t}{\lambda}}$ and we set $v=\frac{v'}{\lambda}$. Here we can identify its vertex set with $\mathbb{Z}_{v}$ and its edges with the pairs $\{x,x+a\}$ such that $x\in \mathbb{Z}_{v}$ and $a$ is not in $\frac{v\lambda}{t}\mathbb{Z}_{v}$, each of which appears $\lambda$ times in the list of edges.
Now we consider the natural projection $\pi$ from $\mathbb{Z}_{v'}$ on $\mathbb{Z}_{\frac{v'}{\lambda}}=\mathbb{Z}_v$.
Clearly $\pi$ is surjective on the vertex set $V\left(^{\lambda}K_{\frac{v'}{t}\times \frac{t}{\lambda}}\right)$. Let $x\in V\left(K_{\frac{v'}{t}\times t}\right)$, then the edges incident with $x$ are mapped into the edges of $E\left(^{\lambda}K_{\frac{v'}{t}\times \frac{t}{\lambda}}\right)$ incident with $\pi(x)$.
Since the graphs $K_{\frac{v'}{t}\times t}$ and $^{\lambda}K_{\frac{v'}{t}\times \frac{t}{\lambda}}$
are both regular with degree $v'-t$,
 those (multi)sets of edges have the same cardinality  and so we can assume that $\pi$ maps the edges incident to $x\in V(\Gamma')$ one-to-one onto the edges incident to $\pi(x)$.
Therefore $\pi$ is a covering map from $K_{\frac{v'}{t}\times t}$ to $^{\lambda} K_{\frac{v'}{t}\times \frac{t}{\lambda}}$.
\end{ex}

Let now consider a graph $\Gamma$, a covering graph $\Gamma'$ of $\G$ and let us denote
by $\pi$ the associated covering map.
Let us suppose that there exists a cellular embedding $\sigma$ of $\Gamma$ in an orientable surface $\Sigma$
associated to the rotation $\rho$ and that there exists a cellular embedding $\sigma'$ of $\Gamma'$ in
an orientable surface $\Sigma'$ associated to the rotation $\rho'$. Because of the face tracing algorithm, see \cite{A},
if we have that $$\pi\circ \rho'((x,y))=\rho \circ \pi((x,y)), \mbox{ for any } \{x,y\}\in E(\Gamma')$$ then $\pi$ maps $\sigma'$-face boundaries into $\sigma$-face boundaries. We remark that $\pi$ does not necessarily preserve the face sizes.
If it does and each edge of $\Gamma$ belongs to two different $\sigma$-faces, $\pi$ induces a homeomorphism between the boundary of a $\sigma'$-face and the corresponding $\sigma$-face boundary. In this case, since the faces are all homeomorphic to an open disc,
$\pi$ can be extended to a map $p: \Sigma'\rightarrow \Sigma$ defined also in the interior of the $\sigma'$-faces of $\Sigma'$ in such a way that, for every $x\in \Sigma$, there exists an open neighborhood $U$ of $x$ such that $p^{-1}(U)$ is a union of disjoint open sets in $\Sigma'$ each of which is mapped homeomorphically onto $U$ by $p$.
Equivalently we can say that $\Sigma'$ is a covering space of $\Sigma$ with respect to the covering map $p$ defined so that the following diagram commutes.
\begin{center}
\begin{tikzcd}
\Gamma' \arrow{r}{\sigma'} \arrow{d}{\pi}
&\Sigma' \arrow{d}{p}\\
\Gamma \arrow{r}{\sigma} &\Sigma
\end{tikzcd}
\end{center}
Now we come back to the embeddings associated to the pair of relative Heffter arrays of Theorem \ref{thm:proiezione}. We can prove the following result.
\begin{thm}
Let $A$ be an $^{\alpha}\H_{t}(m,n;s,k)$ and let  $B$ be the $^{\alpha\lambda} \H_{\frac{t}{\lambda}}(m,n;s,k)$
obtained from $A$  via the projection map
$\pi: \mathbb{Z}_{\frac{2nk}{\alpha}+t}\rightarrow \mathbb{Z}_{\frac{2nk}{\alpha\lambda}+\frac{t}{\lambda}}$.
Let us suppose that $B$ admits simple and compatible orderings $\omega_r$ and $\omega_c$ and denote by
$\sigma:\ ^{\alpha\lambda} K_{(\frac{2nk}{\alpha t}+1)\times \frac{t}{\lambda}}\rightarrow \Sigma$ the associated biembedding.
Then there exists also a biembedding $\sigma'$ of $^{\alpha} K_{(\frac{2nk}{\alpha t}+1)\times t}$ in an orientable surface
$\Sigma'$ that is a covering space for $\Sigma$ with respect to a covering map $p$ such that the following diagram commutes:
\begin{center}
\begin{tikzcd}
^{\alpha} K_{(\frac{2nk}{\alpha t}+1)\times t} \arrow{r}{\sigma'} \arrow{d}{\pi}
&\Sigma' \arrow{d}{p}\\
^{\alpha\lambda} K_{(\frac{2nk}{\alpha t}+1)\times \frac{t}{\lambda}} \arrow{r}{\sigma} &\Sigma
\end{tikzcd}
\end{center}
\end{thm}
\begin{proof}
From the considerations of Example \ref{ex:projection} we have that the map $\pi$ can be seen as a covering map from $^{\alpha} K_{(\frac{2nk}{\alpha t}+1)\times t} $ to $^{\alpha\lambda} K_{(\frac{2nk}{\alpha t}+1)\times \frac{t}{\lambda}}$. We have already observed that if $\omega_r$ and $\omega_c$ are simple and compatible then the orderings $\omega_r'$ and $\omega_c'$ associated to the same permutations $\alpha_r$ and  $\alpha_c$ are simple and compatible also for $A$. Therefore there exists a biembedding $\sigma'$ of $^{\alpha} K_{(\frac{2nk}{\alpha t}+1)\times t}$ in an orientable surface $\Sigma'$. We denote by $\rho$ the rotation associated to $\sigma$ (defined using the permutation $\gamma$ of $skel(B)\times \{1,-1\}$) and by $\rho'$ the rotation associated to $\sigma'$ (defined using the permutation $\gamma'$ of $skel(A)\times \{1,-1\}$). Here we note that the number of filled cells in each row (resp. column) of $A$ is the same as the number of filled cells in each row (resp. column) of $B$. This means that $\pi$ conserves the face sizes. Since $\sigma$ and $\sigma'$ are face $2$-colorable embeddings, we also have that each edge belongs to two different faces.
Therefore, due to the previous discussion, in order to obtain the thesis, it suffices to prove that
$$\pi\circ \rho'((x,y))=\rho \circ \pi((x,y)), \mbox{ for any } \{x,y\}\in E\left(^{\alpha} K_{(\frac{2nk}{\alpha t}+1)\times t} \right).$$
Following the proof of Theorem \ref{thm:biembedding}, each edge of $^{\alpha} K_{(\frac{2nk}{\alpha t}+1)\times t} $
can be written in the form $(x,x+\varepsilon \cdot a_{i,j})$ where $\varepsilon=\pm 1$.
Now given $(x,x+\varepsilon \cdot a_{i,j})\in D\left(^{\alpha} K_{(\frac{2nk}{\alpha t}+1)\times t} \right)$, we have:
\begin{eqnarray}
\nonumber
  \pi\circ \rho'((x,x+\varepsilon \cdot a_{i,j})) &=& \pi((x,x+\gamma_2'(i,j,\varepsilon)\cdot a_{\gamma_1'(i,j,\epsilon)})) \\
\nonumber  \  &=& (\pi(x),\pi(x)+\pi(\gamma_2'(i,j,\epsilon)\cdot a_{\gamma_1'(i,j,\varepsilon)})).
\end{eqnarray}
Now, from the definition of $B$, it follows that $\pi(a_{i,j})=b_{i,j}$ and, since $skel(A)=skel(B)$ and both
$\gamma$ and $\gamma'$ are defined using the same permutations $\alpha_r$ and $\alpha_c$, we have that $\gamma=\gamma'$ that is
$\gamma'_1(i,j,\varepsilon)=\gamma_1(i,j,\varepsilon)$ and $\gamma_2'(i,j,\varepsilon)=\gamma_2(i,j,\varepsilon)$.
Therefore:
$$(\pi(x),\pi(x)+\pi(\gamma_2'(i,j,\epsilon)\cdot a_{\gamma_1'(i,j,\varepsilon)}))=
(\pi(x),\pi(x)+\gamma_2(i,j,\varepsilon)\cdot b_{\gamma_1(i,j,\varepsilon)}).$$
The thesis follows because the right hand side of the previous equality is exactly
$\rho((\pi(x),\pi(x)+\varepsilon \cdot b_{i,j}))=\rho(\pi((x,x+\varepsilon \cdot a_{i,j}))).$
\end{proof}
\section*{Acknowledgements}
The authors were partially supported by INdAM--GNSAGA.

\end{document}